\documentclass[10pt, a4paper, reqno]{amsart}
\usepackage{amsthm,amssymb,amsmath}
\usepackage{graphicx}
\usepackage{enumerate}
\usepackage{hyperref}
\usepackage[latin1]{inputenc}

\theoremstyle{definition}
\newtheorem{definition}{Definition}
\newtheorem{example}[definition]{Example}
\theoremstyle{remark}

\theoremstyle{plain}
\newtheorem{lemma}[definition]{Lemma}

\newtheorem{theorem}[definition]{Theorem}
\newtheorem{cor}[definition]{Corollary}

\newcommand{\set}[1]{\left\{{#1}\right\}}
\newcommand{\vek}[1]{\boldsymbol{#1}}

\newcommand\setsuchas[2]{\left\{\,{#1}\,\vrule\,{#2}\,\right\}}
\newcommand\setsuchasc[2]{\left\{\,{#1}:\,{#2}\,\right\}}

\newcommand{\Z}{{\mathbb{Z}}}
\newcommand{\C}{{\mathbb{C}}}
\newcommand{\R}{{\mathbb{R}}}

\newcommand{\tdeg}[1]{\lvert {#1} \rvert}

\newcommand{\divides}[2]{{#1 \left\lvert {#2} \right.}}
\newcommand{\dividesnot}[2]{{#1 \not \vert #2}}

\newcommand{\degseq}{\boldsymbol{d}}
\newcommand{\sfp}{\operatorname{sfp}}

\begin{document}
\title{Laplacians on shifted multicomplexes}
\author{Jan Snellman}
\address{Department of Mathematics, Linköping University\\
SE-58183 Linköping, Sweden}
\email{Jan.Snellman@liu.se}

\begin{abstract}
  We define the Laplacian operator on finite multicomplexes and give a
  formula for its spectra in the case of shifted multicomplexes.
\end{abstract}
\maketitle

\hyphenation{multi-comp-lexes}
\sloppy

\begin{section}{Introduction}
  The Laplacian of an undirected graph is a square matrix, whose
  eigenvalues yield important information. We can regard graphs as
  one-dimensional simplicial complexes, and ask whether there is a
  generalisation of the Laplacian operator to simplicial complexes. It
  turns out that there is, and that is useful for calculating real
  Betti numbers \cite{Friedman:LaplacianBetti}. 

  Duval and Reiner \cite{DuvalReiner:ShiftedLaplacian} have studied
  Laplacians of a special class of simplicial complexes, the so called
  \emph{shifted} simplicial complexes. They show that such Laplacians
  have integral spectra, computable by a  simple combinatorial
  formula. 

  The present author  \cite{Snellman:Dirichlet} studied a family of monomial
  algebras, corresponding to truncations of the ring of arithmetical
  functions with Dirichlet convolution. The defining monomial ideals
  turn out to be \emph{strongly stable}, which makes it possible to
  use the so-called Eliahou-Kervaire resolution \cite{Eliahou:MinRes}
  to calculate the Betti numbers of the ideals.

  Simplicial complexes correspond to square-free monomial ideals, or
  differently put, to monomial ideals in the exterior algebra. Shifted
  simplicial complexes correspond to strongly stable square-free
  ideals. Monomial ideals, on the other hand, correspond to so-called
  \emph{multicomplexes}. One is  naturally led to the question: is there a
  way to define Laplacian operators on finite multicomplexes, and if
  so, is there a simple formula for their spectra in the case of
  multicomplexes corresponding to strongly stable monomial ideals?

  The Laplacian operator on simplicial complexes is defined as 
  \begin{equation}\label{eq:simplap}
    L_d' = \partial_{d+1} \partial_{d+1}^*,
  \end{equation}
  where \(\partial\) is the simplicial boundary operator, and
  \(\partial^*\) is its dual. Björner and Vre\'cica
  \cite{BjornerMulticomplex} defined a boundary operator on
  multicomplexes which generalises the boundary operator on simplicial
  complexes. Using their definition, we can define Laplacians on
  multicomplexes in analogy with \eqref{eq:simplap}. It turns out that
  the formulas of Duval and Reiner for the spectra of Laplacians of
  shifted simplicial complexes generalise neatly to the case of
  shifted multicomplexes, that is to say, to strongly stable, artinian
  monomial ideals. We finally arrive at a formula for the spectrum of
  the Laplacian of the defining ideals of the truncation algebras
  studied by the author in \cite{Snellman:Dirichlet}.

\end{section}

\begin{section}{Laplacians of multicomplexes}
  Let \(X=\set{x_1,x_2,\dots,x_n}\) be a finite set, which is 
  totally ordered, so that 
  \begin{equation}
    \label{eq:tota}
    x_1 < x_2 < \cdots < x_n
  \end{equation}

  By a
  \emph{multicomplex} on \(X\) 
  we mean a finite subset \(M \subset [X]\) of the free abelian monoid
  on 
  \(X\) which is closed under taking divisors, i.e.
  \begin{equation}
    \label{eq:1}
    \divides{m}{t} \text{ and } t \in M \quad \implies \quad m \in M
  \end{equation}
  We define the  total degree of a monomial in \([X]\) by
  \begin{equation}
    \label{eq:2}
    \tdeg{x_1^{a_1} \cdots x_n^{a_n}} = \sum_{i=1}^n a_i
  \end{equation}
  and let \([X]_\ell\) and \(M_\ell\) denote the subset of monomials in
  \([X]\) and \(M\) of total degree \(\ell\).

  \begin{subsection}{The homology theory of Björner and Vre\'cica}
  Björner and Vre\'cica \cite{BjornerMulticomplex} defined a boundary
  operator \(\partial\) on \(\Z[X]\) by
  \begin{equation}
    \label{eq:multidifferential}
    \begin{split}
      \partial_d: \Z[X]_d & \mapsto \Z[X]_{d-1} \\
      \partial_d(x_{i_0}^{\alpha_0} \cdots x_{i_k}^{\alpha_k}) & =
      \sum_{j=0}^k (-1)^{\alpha_0 + \dots \alpha_{j-1}} \cdot r_j \cdot
        x_{i_0}^{\alpha_0} \cdots x_{i_j}^{\alpha_j-1} \cdots
        x_{i_k}^{\alpha_k}, 
    \end{split}
  \end{equation}
  where \(r_j=0\) if \(\alpha_j\) is even, and \(r_j=1\) if
  \(\alpha_j\) is odd. 
  The map \(\partial\) satisfies \(\partial^2=0\), and since \(M\) is
  closed under taking divisors, \(\partial(\Z M) \subset \Z M\), so
  \(\partial\) is a boundary on \(\Z M\).

  Every monomial \(m \in M\) can be uniquely written
  \begin{equation}
    \label{eq:4}
    m = p^2q
  \end{equation}
  where \(q\) is square-free. We define
  \begin{equation}
    \label{eq:5}
    M^{(p^2)} = \setsuchas{q \in M}{p^2q \in M, \, q \text{
        square-free}} 
  \end{equation}
  Then \(M^{(p^2)}\) is a simplicial complex, and  hence the graded
  \(\Z\)-module \(\Z M\) decomposes as
  \begin{equation}
    \label{eq:dec1}
    \Z M = \bigoplus_{p^2 \in M}  p^2 \Z M^{(p^2)}
  \end{equation}
  or, taking the grading into account,
  \begin{equation}
    \label{eq:dec2}
     \Z M_{\ell} = 
     \bigoplus_{p^2 \in M, \, 2\tdeg{p} \le \ell} 
     p^2 \Z M_{\ell-2\tdeg{p}}^{(p^2)}
  \end{equation}
  where
  \begin{equation}
    \label{eq:9}
    M_{\nu}^{(p^2)} = \setsuchas{q \in M}{p^2q \in M, \, q \text{
        square-free}, \, \tdeg{q}=\nu} 
  \end{equation}
  
  To simplify matters, we have shifted the ordinary grading on
  simplicial 
  complexes so that the empty simplex lies in dimension zero rather
  than in dimension \(-1\). Defining the \(\boldsymbol{f}\)-vector of
  a multicomplex as the vector \(\boldsymbol{f}(M)=(f_0,f_1,\dots)\),
  where \(f_i = \# M_i\), we get that
  \begin{equation}
    \label{eq:10}
    \boldsymbol{f}(M) = \sum_{p^2 \in M}
    S^{2\tdeg{p}}\boldsymbol{f}(M^{p^2}),  
  \end{equation}
  where \(\boldsymbol{f}(M^{p^2})\) is the re-indexed
  \(\boldsymbol{f}\)-vector of the simplicial complex \(M^{p^2}\), and
  \(S: \Z^\infty \to \Z^\infty\) is the shift map
  \begin{equation}
    \label{eq:11}
    S((f_0,f_1,f_2,\dots)) = (0,f_0,f_1,f_2,\dots)
  \end{equation}

  Since
  \begin{equation}
    \label{eq:6}
    \partial(p^2q) = p^2\partial(q)
  \end{equation}
  the boundary operator on \(M\) restricts to the ordinary simplicial
  boundary operator on \(M^{(p^2)}\). Hence, the decomposition
  \eqref{eq:dec1} is a decomposition of differential graded
  \(\Z\)-modules. It follows that the homology of \(M\) splits as
  \begin{equation}
    \label{eq:homologysplit}
    H_{\ell}(M) \simeq \bigoplus_{p^2 \in M, \, 2\tdeg{p} \le \ell} 
    H_{\ell-2\tdeg{p}}(M^{(p^2)})
  \end{equation}

  \end{subsection}

  \begin{subsection}{Definition of Laplacians on multicomplexes}

  We now proceed to define the \emph{Laplacian} of a multicomplex. It
  is convenient to work with real coefficients, so we extend scalars
  to \(\R\) and regard \(\partial_d\) as (the restriction to \(\R M\) of)
  a map 
  \begin{math}
    \partial_d: \, \R[X]_d  \to \R[X]_{d-1}.
  \end{math}
  There is a dual map
  \begin{equation}
    \label{eq:dualmap}
    \partial_{d+1}^*: \, \textrm{Hom}_\R(M_d,\R)  \to
    \textrm{Hom}_\R(M_{d+1},\R) 
  \end{equation}
  We identify \(\textrm{Hom}_\R(M_d,\R) \simeq M_d\), so that
  \begin{displaymath}
    \partial_{d+1}^*: \, \R M_d  \to \R M_{d+1}
  \end{displaymath}
  With respect to the natural basis, the matrix of \(\partial_{d+1}^*\)
  is just the 
  transpose of the matrix of \(\partial_d\). 
  More explicitly, if \(m=x_1^{a_1}\cdots x_n^{a_n}\) then
  \begin{equation}
    \label{eq:dstar}
    \begin{split}
    \partial_{d+1}^*(m) &= \sum_{j=1}^n (-1)^{a_1 + \dots + a_{j-1}}
    s_j x_j m, \\
    s_j &=
    \begin{cases}
      1 & \text{ if } a_j \text{ is even and } x_j m \in M \\
      0 & \text{ otherwise}
    \end{cases}
    \end{split}
  \end{equation}

  \begin{definition}\label{def:laplacian}
  \begin{equation}
    \label{eq:laplacian}
    \begin{split}
      L_d' &= \partial_{d+1} \partial_{d+1}^* \\
      L_d'' &= \partial_{d}^* \partial_{d} \\
      L_d &= L_d' + L_d''
    \end{split}
  \end{equation}
  The map \(L_d\) is the \emph{Laplacian} of the multicomplex
  \(M\).    
  \end{definition}
  While \(L_d\) is the object which interest us, it will for technical
  reasons be more convenient to study the map \(L_d'\). As we
  shall see, the spectrum of \(L_d'\) determines that of \(L_d\), and
  vice versa, so no information is lost by this change of focus.

  Definition \eqref{def:laplacian} is, \emph{mutatis mutandis}, the
  definition of Laplacians for 
  simplicial complexes used by Duval and Reiner \cite{Duval:RecursionLaplacian,
    DuvalReiner:ShiftedLaplacian}. We will regard simplicial complexes
  as special cases of multicomplexes, namely as multicomplexes
  consisting only of square-free monomials. The definitions we make for
  multicomplexes (e.g., of the spectra of their Laplacians) will then
  specialise to the (well studied) case of simplicial complexes.  
  \end{subsection}

 \end{section}

\begin{section}{The spectrum of the Laplacian}\label{sec:SpectrumLaplacian}
  From \eqref{eq:6} it follows that   
  \begin{equation}
    \label{eq:lap6}
    \begin{split}
      L'(p^2q) & = p^2L'(q) \\
      L''(p^2q) & = p^2L''(q) \\
      L(p^2q) & = p^2L(q) 
    \end{split}
  \end{equation}
  Following Duval \cite{Duval:RecursionLaplacian,
    DuvalReiner:ShiftedLaplacian} we define the spectra 
  \(\boldsymbol{s}_i'\),   \(\boldsymbol{s}_i''\),
  \(\boldsymbol{s}_i^{tot}\),
  of the selfadjoint, nonnegative definite operators
  \(L_i'\),  \(L_i''\), \(L_i\)
  to be the multiset  of their (real and
  and nonnegative) eigenvalues. We will identify such a multiset with its
  weakly decreasing rearrangement, which is a partition, and we will,
  unless otherwise stated, identify such partitions that differ only
  in the number of zero parts.

  The spectrum of the Laplacian of \(M\) is the multiset sum of the
  spectra of the constituent simplicial complexes \(M^{p^2}\).

  \begin{lemma}\label{lemma:dirsum}
  \begin{equation}
    \label{eq:spectrasplit}
    \begin{split}
      \boldsymbol{s}_i^{'}(M,\partial) &= \sum_{p^2 \in M, \, 2
        \deg(p) \le i}
      \boldsymbol{s}_{i-2\deg(p)}^{'}(M^{p^2},\partial) \\
      \boldsymbol{s}_i^{''}(M,\partial) &= \sum_{p^2 \in M, \, 2
        \deg(p) \le i} 
      \boldsymbol{s}_{i-2\deg(p)}^{''}(M^{p^2},\partial) \\
      \boldsymbol{s}_i^{tot}(M,\partial) &= \sum_{p^2 \in M, \, 2
        \deg(p) \le i}
      \boldsymbol{s}_{i-2\deg(p)}^{tot}(M^{p^2},\partial) 
    \end{split}
  \end{equation}
  \end{lemma}
  
  Since the spectra of the Laplacians of the constituent simplicial
  complexes are   invariant under reordering of the vertices
  \cite[Remark 3.2]{DuvalReiner:ShiftedLaplacian}, the same
  holds true for  multicomplex Laplacian spectra. Furthermore, the
  relations
  \begin{equation}
    \label{eq:relsLap}
    \begin{split}
      \boldsymbol{s}_i^{''}(M^{p^2},\partial) &= 
      \boldsymbol{s}_{i-1}^{'}(M^{p^2},\partial) \\
      \boldsymbol{s}_i^{tot}(M^{p^2},\partial) &= 
      \boldsymbol{s}_{i}^{'}(M^{p^2},\partial) \cup
      \boldsymbol{s}_{i}^{''}(M^{p^2},\partial)  \\
      \boldsymbol{s}_i^{'}(M^{p^2},\partial) &= 
      \boldsymbol{s}_{i}^{tot}(M^{p^2},\partial) -
      \boldsymbol{s}_{i}^{''}(M^{p^2},\partial)
    \end{split}
  \end{equation}
  which holds for all simplicial complexes \(M^{p^2}\) according to
  \cite{DuvalReiner:ShiftedLaplacian}, translates using
  Lemma~\ref{lemma:dirsum} to
  \begin{lemma}
  \begin{equation}
    \label{eq:relsLap}
    \begin{split}
      \boldsymbol{s}_i^{''}(M,\partial) &= 
      \boldsymbol{s}_{i-1}^{'}(M,\partial) \\
      \boldsymbol{s}_i^{tot}(M,\partial) &= 
      \boldsymbol{s}_{i}^{'}(M,\partial) \cup
      \boldsymbol{s}_{i}^{''}(M,\partial)  \\
      \boldsymbol{s}_i^{'}(M,\partial) &= 
      \boldsymbol{s}_{i}^{tot}(M,\partial) -
      \boldsymbol{s}_{i}^{''}(M,\partial)
    \end{split}
  \end{equation}
  Here, equality means equality as multisets, except that the
  multiplicity of the element zero may differ; \(\cup\) means union of
  multisets, i.e. all occurring multiplicities are added, and \(-\)
  denotes multiset difference.
  \end{lemma}

  The lemma shows that the spectra of the Laplacians
  \(\setsuchas{L_i}{i \ge 0}\) is
  completely determined by the eigenvalues of the operators
  \(\setsuchas{L_i'}{i \ge 0}\). 
\end{section}

\begin{section}{Shifted multicomplexes and their Laplacian spectra}
  We will give a combinatorial formula for the spectrum of \(L'\) for
  the special case of shifted multicomplexes.
  
  \begin{definition}
    A subcomplex \(N \subseteq M\) is \emph{shifted} (relative its support) if
    \begin{equation}
      \label{eq:shifted}
      x_jm \in N, \quad i < j, \quad x_i \in N \qquad
      \implies \qquad x_i m \in N
    \end{equation}
  \end{definition}
  
  \begin{lemma}
    If \(M\) is shifted, then so are all \(M_{p^2}\), with the induced
    total ordering on the vertices in their supports.
  \end{lemma}
  \begin{proof}
    Let \(x_j q \in M_{p^2}\), so that \(x_j q p^2 \in M\) and \(x_j
    q\) and \(p\) have disjoint support. Now, if \(i < j\) and \(x_i
    \in M_{p^2}\) then \(x_i \in M\) and \(\dividesnot{x_i}{p}\).
    Hence, since \(M\) is
    shifted, \(x_i q p^2 \in M\). It follows that \(x_i q \in
    M_{p^2}\). 
  \end{proof}
  
  \begin{definition}
    Let \(N \subseteq M\) be a multicomplex.
    The \emph{degree-sequence} \(\degseq_k\) is the sequence
    \begin{equation}
      \label{eq:degseq}
      \degseq_k(N) = (d_1,d_2,d_3,\dots,d_n)
    \end{equation}
    where \(d_j\) denotes the number of monomials in \(N_k\)
    that are divisible by \(x_j\).
  \end{definition}
  
  \begin{lemma}
    If \(N\) is shifted then \(\degseq_k(N)\) is weakly decreasing,
    i.e. a partition.
  \end{lemma}
  
  \begin{theorem}[Reiner and Duval]
    If \(N\) is a shifted simplicial complex, then
    \begin{equation}
      \label{eq:spectrumlap}
      \boldsymbol{s}_k^{'} = \degseq_k^T(N),
    \end{equation}
    where \(\degseq_k^T(N)\) is the conjugate partition to
    \(\degseq_k(N)\).

    In particular, the eigenvalues of \(L_d'\) are non-negative
    integers. 
  \end{theorem}

  We conclude:
  \begin{theorem}
    Suppose that \(M\) is shifted. Then
    \begin{equation}
      \label{eq:masterthm}
      \boldsymbol{s}_k^{'} = \sum_{p^2 \in M, \, 2
        \deg(p) \le k} \degseq_k^T(M_{p^2}) 
    \end{equation}
    In particular, the eigenvalues of \(L_d'\) are non-negative
    integers. 

    Equivalently: let \(\underline{b}\) be 1 if \(b\) is odd, and zero
    otherwise, and let 
    \begin{equation}
      \label{eq:mod}
      \underline{(\alpha_1,\dots,\alpha_n)}  = 
      (\underline{\alpha_1} ,\dots,\underline{\alpha_n} )
    \end{equation}
    Then
    \begin{equation}
      \label{eq:mthm}
      {\boldsymbol{s}_k^{'}}^T= \sum_{\vek{x}^{\vek{\alpha}} \in M_k}
      \underline{\vek{\alpha}}
    \end{equation}
  \end{theorem}

  \begin{example}
    Let \(M_3=[x_1,x_2,x_3]_3\). Then the matrix of \(d_3\), with
    respect to the basis of monomials 
    of degree three and two ordered lexicographically, is
    \begin{displaymath}
      P=
    \begin{array}{c|cccccccccc}
        & x_1^3 & x_1^2x_2 & x_1^2x_3 & x_1x_2^2 &
        x_1x_2x_3 &
        x_1x_3^2 & x_2^3 & x_2^2 x_3  & x_2x_3^2 &
        x_3^3 \\ \hline
        x_1^2   & 1 & -1 & -1 & 0 & 0  & 0 & 0 & 0 & 0 & 0       \\
        x_1x_2  & 0 &  0 &  0 & 0 & 1 & 0 & 0 & 0 & 0 & 0       \\
        x_1x_3  & 0 &  0 &  0 & 0 & -1  & 0 & 0 & 0 & 0 & 0      \\
        x_2^2   & 0 &  0 &  0 & 1 & 0  & 0 & 1 & -1& 0 & 0      \\
        x_2x_3  & 0 &  0 &  0 & 0 & 1 & 0 & 0 & 0 & 0 & 0     \\
        x_3^2   & 0 &  0 &  0 & 0 & 0  & 1 & 0 & 0 & 1 & 1
    \end{array}
    \end{displaymath}
    The symmetric matrix \(PP^*\) has eigenvalues  \(3,3,3,3,0,0\).
    We have that
    \begin{displaymath}
      \begin{split}
      (3,3,3,3)^T &= (4,4,4) \\
      &= (1,0,0) + (0,1,0) + (0,0,1) + (1,0,0) + (1,1,1) + \\
      & \qquad   (1,0,0) + (0,1,0) + (0,0,1) + (0,1,0) + (0,0,1).        
      \end{split}
    \end{displaymath}

    Let now
    \[
    M_3=\set{x_1^3,x_1^2x_2,x_1x_2^2,x_2^3,x_2^2x_3,x_1x_2x_3,x_1^2x_3}
      \subset [x_1,x_2,x_3]_3,\] 
      as in the left part of Table~\ref{tab:multiex}.
      
      \begin{table}[b]
        \centering
        \begin{tabular}{cc}
    \setlength{\unitlength}{0.8cm}
    \begin{picture}(6,5)(0,-0.8)
      \multiput(0,0)(1,0){4}{\circle{0.3}}
      \multiput(0.5,1)(1,0){3}{\circle{0.3}}
      \multiput(1,2)(1,0){2}{\circle{0.3}}
      \multiput(1.5,3)(1,0){1}{\circle{0.3}}

      \put(-0.9,-0.1){\(x_2^3\)}
      \put(3.5,-0.1){\(x_3^3\)}
      \put(1.1,3.5){\(x_1^3\)}

      \multiput(0,0)(0.5,1){4}{\circle*{0.3}}
      \multiput(1,0)(0.5,1){3}{\circle*{0.3}}
    \end{picture}
&
    \setlength{\unitlength}{0.8cm}
    \begin{picture}(6,5)(0,-0.8)
      \multiput(0,0)(1,0){4}{\circle{0.3}}
      \multiput(0.5,1)(1,0){3}{\circle{0.3}}
      \multiput(1,2)(1,0){2}{\circle{0.3}}
      \multiput(1.5,3)(1,0){1}{\circle{0.3}}

      \put(-0.9,-0.1){\(x_2^3\)}
      \put(3.5,-0.1){\(x_3^3\)}
      \put(1.1,3.5){\(x_1^3\)}

      \multiput(0,0)(0.5,1){4}{\circle*{0.3}}
      \multiput(1.5,1)(0.5,1){2}{\circle*{0.3}}
    \end{picture}
      \end{tabular}
      \caption{Multicomplexes on three variables, degree 2}
      \label{tab:multiex}
    \end{table}

      Now, the matrix of \(d_3\), with
    respect to the basis of monomials 
    of degree three and two ordered lexicographically, is
    \begin{displaymath}
      P=
    \begin{array}{c|cccccccccc}
        & x_1^3 & x_1^2x_2 & x_1^2x_3 & x_1x_2^2 &
        x_1x_2x_3 &
        x_1x_3^2 & x_2^3 & x_2^2 x_3  & x_2x_3^2 &
        x_3^3 \\ \hline
        x_1^2   & 1 & -1 & -1 & 0 & 0 & x & 0 & 0 & x & x       \\
        x_1x_2  & 0 &  0 &  0 & 0 & 1 & x & 0 & 0 & x & x       \\
        x_1x_3  & 0 &  0 &  0 & 0 & -1& x & 0 & 0 & x & x      \\
        x_2^2   & 0 &  0 &  0 & 1 & 0 & x & 1 & -1& x & x      \\
        x_2x_3  & 0 &  0 &  0 & 0 & 1 & x & 0 & 0 & x & x     \\
        x_3^2   & 0 &  0 &  0 & 0 & 0 & x & 0 & 0 & x & x
    \end{array}
    \end{displaymath}
    and \(PP^*\) has eigenvalues  \(3,3,3,0,0,0\).
    We have that
    \begin{displaymath}
      \begin{split}
      (3,3,3)^T &= (3,3,3) \\
      &= (1,0,0) + (0,1,0) + (0,0,1) + (1,0,0) + (1,1,1) + \\
      & \qquad   \quad + (0,1,0) + (0,0,1).
      \end{split}
    \end{displaymath}

  Finally, let
    \[
    M_3=\set{x_1^3,x_1^2x_2,x_1x_2^2,x_2^3,x_1x_2x_3,x_1^2x_3}
      \subset [x_1,x_2,x_3]_3,\] 
      as in the right part of Table~\ref{tab:multiex}.

   Then the matrix of \(d_3\), with    respect to the basis of monomials
    of degree three and two ordered lexicographically, is
    \begin{displaymath}
      P =
    \begin{array}{c|cccccccccc}
        & x_1^3 & x_1^2x_2 & x_1^2x_3 & x_1x_2^2 &
        x_1x_2x_3 &
        x_1x_3^2 & x_2^3 & x_2^2 x_3  & x_2x_3^2 &
        x_3^3 \\ \hline
        x_1^2   & 1 & -1 & -1 & 0 &0  & x & 0 & x & x & x       \\
        x_1x_2  & 0 &  0 &  0 & 0 &1  & x & 0 & x & x & x       \\
        x_1x_3  & 0 &  0 &  0 & 0 &-1 & x & 0 & x & x & x      \\
        x_2^2   & 0 &  0 &  0 & 1 &0  & x & 1 & x & x & x      \\
        x_2x_3  & 0 &  0 &  0 & 0 &1  & x & 0 & x & x & x     \\
        x_3^2   & 0 &  0 &  0 & 0 &0  & x & 0 & x & x & x
    \end{array}
    \end{displaymath}
     \(PP^*\) now has eigenvalues  \(3,3,2,0,0,0\).
    We have that
    \begin{displaymath}
      \begin{split}
      (3,3,2)^T &= (3,3,2) \\
      &= (1,0,0) + (0,1,0) + (0,0,1) + (1,0,0) + (1,1,1) + 
       (0,1,0).
      \end{split}
    \end{displaymath}
  \end{example}

\end{section}

\begin{section}{Strongly stable monomial ideals} 

  Let \(I \subset \R[x_1,\dots,x_n]\) be an artinian monomial ideal,
  i.e., a monomial ideal such that the quotient ring \(S=R/I\) is
  artinian. Let \(M\) be the set of monomials in \([x_1,\dots,x_n]\)
  not in \(I\). Then \(M\) is a finite multicomplex on \([x_1,\dots,x_n]\).
  
  Conversely, if \(M\) is a finite multicomplex on
  \([x_1,\dots,x_n]\), then one can form the monomial ideal \(I\)
  generated by the monomials not in \(I\), and \(M\) will be an \(\R\)-basis
  for the artinian ring   \(S=R/I\). We denote this ring by \(\R M\)
  and call it the multicomplex ring of \(M\). This is coherent with
  our previous use of \(\R M\); we have merely introduced a
  multiplication on this \(\R\)-vector space.

  By means of the correspondence \(M \leftrightarrow I\) between
  finite multicomplexes and artinian monomial ideals, we can now define
  Laplacians and their spectra for artinian monomial ideals. For a
  special class of monomial ideals, which is important because of its
  connection to so-called \emph{generic initial ideals} (see for
  instance \cite{Green:gin,Ebud:View,Bigatti:Betti}), we can use the results of the
  previous section to calculate these spectra.
  
  A monomial ideal \(I\) is said to be
  \emph{strongly stable} or \emph{Borel-fixed} \cite{Green:gin,
    Ebud:View} w.r.t. a total order of \(X=\set{x_1,\dots,x_n}\) 
  if
  \begin{equation}
    \label{eq:shiftedIdeal}
    m \in I \text{ and } \divides{x_i}{m} \text{ and } x_j < x_i \quad
    \implies \quad \frac{x_j}{x_i} m \in I
  \end{equation}

The minimal graded free resolution of strongly stable ideals are
  given by the so-called Eliahou-Kervaire resolution
  \cite{Eliahou:MinRes}. This means that the graded Betti numbers of
  \(I\), and hence of \(S\), can be easily read of from the minimal
  monomial generators of \(I\).

  \begin{lemma}
    \(I\) is strongly stable if and only if \(M\) is shifted w.r.t. the
    reverse ordering on \(X\).
  \end{lemma}
  
  The above lemma means that the spectra of Laplacians of artinian strongly stable
  monomial ideals can be easily calculated from combinatorial data involving the
  standard monomials in the quotient ring.
  
\end{section}

\begin{section}{Truncations of arithmetical functions with Dirichlet convolution}
  We will apply the results of section~\ref{sec:SpectrumLaplacian} to
  a particular family of shifted multicomplexes that arise naturally
  in the study of arithmetical functions. Recall that an arithmetical
  function is a complex-valued function defined on the positive
  integers. The set \(\Gamma\) of all such functions becomes a complex
  vector space under point-wise addition and multiplication by
  scalars. If we introduce the indicator functions \(e_m\), defined
  for all positive integers \(m\) by
  \begin{equation}
    \label{eq:em}
    e_m(k) =
    \begin{cases}
      1 & k=m \\ 0 & k \neq m
    \end{cases}
  \end{equation}
  then every element \(f \in \Gamma\) can be written as a formal
  linear combination
  \begin{equation}
    \label{eq:flin}
    f = \sum_{m=1}^\infty c_m e_m, \qquad c_m \in \C
  \end{equation}
  With respect to a natural adic topology on \(\Gamma\), the above is
  an absolutely convergent sum.

  There is also a famous convolution product, the so-called
  \emph{Dirichlet convolution}, with respect to which \(\Gamma\)
  becomes an associative, commutative \(\C\)-algebra isomorphic to the
  unrestricted power series ring on countably many variables.
  
  The Dirichlet convolution is defined by
  \begin{equation}
    \label{eq:Dirichlet}
    f*g(m) = \sum_{\divides{k}{m}} f(k)g(m/k)
  \end{equation}
  which is just the \(\C\)-linear and continuous extension of the rule
  \begin{equation}
    \label{eq:rule}
    e_a*e_b = e_{ab}
  \end{equation}
  The isomorphism between \((\Gamma,*)\) and  \(\C[[x_1,x_2,\dots,]]\)
  alluded to above is given by the \(\C\)-linear and continuous
  extension of
  \begin{equation}
    \label{eq:isom}
    e_{p_1^{a_1} \cdots p_r^{a_r}} \mapsto 
    y_1^{a_1} \cdots y_r^{a_r}
  \end{equation}
  where \(p_i\) is the \(i\)'th prime number. Cashwell and Everett
  \cite{Cashwell:FPS} used this isomorphism to show that
  \((\Gamma,*)\) is a unique factorisation domain.

  Now let \(N\) be a positive integer, and let \(\Gamma_N\) denote the
  subset of all arithmetical functions supported on
  \(\set{1,2,\dots,N}\). With the modified multiplication
  \begin{equation}
    \label{eq:dirmod}
    f *_N g (k) =
    \begin{cases}
      f*g(k) & k \le N \\
      0 & k > N
    \end{cases}
  \end{equation}
  \(\Gamma_N\) becomes a retract of \(\Gamma\). Furthermore,
  \(\Gamma_N\) is the multicomplex ring on \(M_N\), the multicomplex
  consisting of \(\set{1,2,\dots, N}\), where \(k=p_1^{a_1} \cdots
  p_r^{a_r}\) is regarded as a multiset on \(\set{p_1,\dots,p_r}\),
  i.e. is identified with the monomial \(x_1^{a_1} \cdots
  x_r^{a_r}\). 

  In \cite{Snellman:Dirichlet} the multicomplex ring \(\Gamma_N\),
  regarded as the artinian monomial algebra \(\C[x_1,\dots,x_n]/I_N\),
  was studied. The point of departure was the fact that \(I_N\) is
  strongly stable w.r.t. the reverse order of the variables, so that
  the Eliahou-Kervaire resolution yields the graded Betti numbers of
  the algebra in terms of certain number-theoretic quantities. We will
  here use the fact that \(M_N\) is a shifted multicomplex to
  determine its Laplacian spectrum.

  \begin{lemma}
    \(M_N\) is a shifted multicomplex on \(y_1,\dots,y_n\), where
    \(n\) is the largest number so that \(p_n \le N\).
  \end{lemma}
  \begin{proof}
    Obvious.
  \end{proof}

  \begin{definition}
    For a positive integer \(m\) with prime factorisation
    \(m=p_1^{a_1} \cdots p_r^{a_r}\) we define
    \(\log(m)=(a_1,a_2,\cdots)\). Conversely, given a finitely
    supported sequence \(\vek{\alpha}=(a_1,a_2,a_3,\dots)\) of
    non-negative integers, we define \(\vek{p}^{\vek{\alpha}} =
    p_1^{a_1} \cdots p_r^{a_r}\). 
  \end{definition}

  With this definition, the \emph{squarefree part} of
  \(\vek{p}^{\vek{\alpha}}\) is
  \[\sfp(\vek{p}^{\vek{\alpha}})=\vek{p}^{\underline{\vek{\alpha}}},\]
  and the total degree
  \(\tdeg{\vek{\alpha}} = \sum_i a_i\) is equal to the number of not
  necessarily distinct prime factors of \(\vek{p}^{\vek{\alpha}}\),
  i.e. to \(\Omega(\vek{p}^{\vek{\alpha}})\) to use the notation in \cite{HW}.
  
  \begin{cor}
    The spectrum \(\boldsymbol{s}_k^{'}\) of the \(k\)'th Laplacian
    \(L_k'\) of the multicomplex \(M_N\) is given by 
    \begin{equation}
      \label{eq:MnLaplSpec}
      {\boldsymbol{s}_k^{'}}^T = 
      \sum_{\substack{1 \le \ell \le N\\
          \Omega(\ell)=k}} \underline{\log(\ell)} 
      = 
      \sum_{\substack{1 \le \ell \le N\\
          \Omega(\ell)=k}} \log(\sfp(\ell)) 
    \end{equation}
    Hence, if we introduce the arithmetical functions \(t_k^i\) and
    \(s_k^i\) by
    \begin{equation}
      \label{eq:sL}
      \begin{split}
        {\boldsymbol{s}_k^{'}}(N)  &= (t_k^1(N),t_k^2(N),\dots)  \\      
        {{\boldsymbol{s}_k^{'}}(N)}^T  &= (s_k^1(N),s_k^2(N),\dots)  
      \end{split}
    \end{equation}
    
    then
    \begin{equation}
      \label{eq:sIs}
      \begin{split}
      t_k^i(N) &= \sum_{\substack{1 < n \le N\\ \Omega(n)=k \\
          \divides{p_i}{\sfp(n)}}} 1   \\
      s_k^j(N) &= \sum_{\setsuchas{i}{t_k^i(N) \ge j}} 1 \\
      &= \# \setsuchasc{\ell}{\# \setsuchasc{1 < n \le
          N}{\Omega(n)=k, \,\divides{p_\ell}{\sfp(n)}} \ge j}
      \end{split} 
    \end{equation}
  \end{cor}

  We have that \(s_1^1(N)\) is the number of primes \(\le N\), and
  that \(s_1^j(N)=0\) for \(j>1\).

  To study \(t_2^i(N)\) and \(s_2^i(N)\), first note that 
  \begin{equation}
    \label{eq:7}
    \sfp(p_ap_b) =
    \begin{cases}
      p_ap_b & a \neq b \\
      1 & a = b
    \end{cases}
  \end{equation}
  It follows that if we define
  \begin{equation}
    \label{eq:Y2}
    \begin{split}
      Y_2(N)  &= (y_{ab}), \qquad y_{ab}=
      \begin{cases}
        1 & a \neq b, \, p_a p_b \le N \\
        0 & \text{ otherwise}
      \end{cases}   \\
      U_2(N) &= (u_{ab}), \qquad u_{ab} =
      \begin{cases}
        y_{ab} & b < a \\
        y_{a,b+1} & b \ge a
      \end{cases}
    \end{split}
  \end{equation}
  then \(t_2^i(N)\) is the \(i\)'th row sum of \(U_2(N)\), 
  and \(s_2^j(N)\) is the \(j\)'th column sum of \(U_2(N)\).
  An example, for \(N=50\), is shown in Table~\ref{tab:Y2.50}.

  \begin{table}[t]
    \centering
    \begin{tabular}{cc}
      \begin{math}
      \begin{bmatrix}
          0&  1&  1&  1&  1&  1&  1&  1&  1 \\
          1&  0&  1&  1&  1&  1&  0&  0&  0 \\
          1&  1&  0&  1&  0&  0&  0&  0&  0 \\
          1&  1&  1&  0&  0&  0&  0&  0&  0 \\
          1&  1&  0&  0&  0&  0&  0&  0&  0 \\
          1&  1&  0&  0&  0&  0&  0&  0&  0 \\
          1&  0&  0&  0&  0&  0&  0&  0&  0 \\
          1&  0&  0&  0&  0&  0&  0&  0&  0 \\
          1&  0&  0&  0&  0&  0&  0&  0&  0  
      \end{bmatrix}
      \end{math}
      &
      \begin{math}
      \begin{bmatrix}
          1&  1&  1&  1&  1&  1&  1&  1&  0 \\
          1&  1&  1&  1&  1&  0&  0&  0& 0 \\
          1&  1&   1&  0&  0&  0&  0&  0&0 \\
          1&  1&  1&  0&  0&  0&  0&  0&  0 \\
          1&  1&  0&  0&  0&  0&  0&  0&  0 \\
          1&  1&  0&  0&  0&  0&  0&  0&  0 \\
          1&  0&  0&  0&  0&  0&  0&  0&  0 \\
          1&  0&  0&  0&  0&  0&  0&  0&  0 \\
          1&  0&  0&  0&  0&  0&  0&  0&  0 
      \end{bmatrix}
      \end{math}
    \end{tabular}
    \caption{\(Y_2(50)\) and \(U_2(50)\). \(\mathbf{t}_2=(8,5,3,3,2,2,1,1,1)\),
      \(\mathbf{s}_2=(9,6,4,2,2,1,1,1)\).}
    \label{tab:Y2.50}
  \end{table}

  Since \(p_n \sim n \log n\), 
  we see that for fixed \(i\) and large \(N\),
  \begin{equation}
    \label{eq:as1}
    s_2^{i}(N) \sim \frac{N/p_i}{\mathcal{W}(N/p_i)}
  \end{equation}
  where \(\mathcal{W}\) is the Lambert W-function. If we allow also
  \(i\) to tend to infinity, but much slower than \(N\), then we get
  that
  \begin{equation}
    \label{eq:as2}
    s_2^{i}(N) \sim \frac{N/(i \log i)}{\mathcal{W}(N/(i \log i))}
  \end{equation}
\end{section}

\bibliographystyle{plain}
\bibliography{journals,snellman,CombinatoricsArticles,CommutativeAlgebraArticles,NumberTheoryArticles}

\def\cprime{$'$}
\begin{thebibliography}{10}

\bibitem{Bigatti:Betti}
A.~M. Bigatti.
\newblock Upper bounds for the {B}etti numbers of a given {H}ilbert function.
\newblock {\em Communications in {A}lgebra}, 21:2317--2334, 1993.

\bibitem{BjornerMulticomplex}
Anders Bj{\"o}rner and Sini{\v{s}}a Vre{\'c}ica.
\newblock On {$f$}-vectors and {B}etti numbers of multicomplexes.
\newblock {\em Combinatorica}, 17(1):53--65, 1997.

\bibitem{Cashwell:FPS}
E.~D. Cashwell and C.~J. Everett.
\newblock Formal power series.
\newblock {\em Pacific {J}ournal of {M}athematics}, 13:45--64, 1963.

\bibitem{Duval:RecursionLaplacian}
Art~M. Duval.
\newblock A common recursion for laplacians of matroids and shifted simplicial
  complexes.
\newblock {\em Documenta Math.}, 10:583--618, 2005.

\bibitem{DuvalReiner:ShiftedLaplacian}
Art~M. Duval and Victor Reiner.
\newblock Shifted simplicial complexes are {L}aplacian integral.
\newblock {\em Trans. Amer. Math. Soc.}, 354(11):4313--4344 (electronic), 2002.

\bibitem{Ebud:View}
David Eisenbud.
\newblock {\em Commutative {A}lgebra with a {V}iew {T}oward {A}lgebraic
  {G}eometry}, volume 150 of {\em Graduate {T}exts in {M}athematics}.
\newblock Springer {V}erlag, 1995.

\bibitem{Eliahou:MinRes}
S.~Eliahou and M.~Kervaire.
\newblock Minimal resolutions of some monomial ideals.
\newblock {\em Journal of {A}lgebra}, 129:1--25, 1990.

\bibitem{Friedman:LaplacianBetti}
J.~Friedman.
\newblock Computing {B}etti numbers via combinatorial {L}aplacians.
\newblock {\em Algorithmica}, 21(4):331--346, 1998.

\bibitem{Green:gin}
Mark~L. Green.
\newblock Generic initial ideals.
\newblock In {\em Proceedings of the {S}ummer {S}chool on {C}ommutative
  {A}lgebra}, volume~2, pages 16--85, CRM, Barcelona, 1996.

\bibitem{HW}
G.~H. Hardy and E.~M. Wright.
\newblock {\em An introduction to the theory of numbers}.
\newblock The Clarendon Press Oxford University Press, New York, fifth edition,
  1979.

\bibitem{Snellman:Dirichlet}
Jan Snellman.
\newblock Truncations of the ring of number-theoretic functions.
\newblock {\em Homology Homotopy Appl.}, 2:17--27 (electronic), 2000.

\end{thebibliography}

\end{document}